\documentclass[12pt]{article}

\usepackage{color,latexsym,amsfonts,amssymb}
\usepackage{epsfig}
\usepackage{amsthm}
\usepackage{amsmath}

\date{}
\begin{document}
 \def \be {\begin{equation}}
 \def \ee {\end{equation}}
 \def \beq {\begin{eqnarray}}
 \def \eeq {\end{eqnarray}}
 \def \beqs {\begin{eqnarray*}}
 \def \eeqs {\end{eqnarray*}}
 \def \bt {\begin{theorem}}
 \def \et {\end{theorem}}
 \def \bl {\begin{lemma}}
 \def \el {\end{lemma}}
 \def \bc {\begin{corollary}}
 \def \ec {\end{corollary}}
 \def \bd {\begin{definition}}
 \def \ed {\end{definition}}
 \def \qed {\hfill $\spadesuit $}
 \def \bp {\begin{proposition}}
 \def \ep {\end{proposition}}
 \def \bcon {\begin{conjecture}}
 \def \econ {\end{conjecture}}
 \def \bq {\begin{question}}
 \def \eq {\end{question}}
\newtheorem{theorem}{Theorem}
\newtheorem{lemma}{Lemma}
\newtheorem{corollary}{Corollary}
\newtheorem{definition}{Definition}
\newtheorem{proposition}{Proposition}
\newtheorem{conjecture}{Conjecture}
\newtheorem{question}{Question}

\title{Uniformly bounded components of normality}

\author{Xiaoling WANG\footnote{Nanjing University of Finance and Economics}
\footnote{This Project was granted financial support from
Postdoctoral Foundation of Shantou University} \  and \ Wang
ZHOU\footnote{National University of Singapore}
 \footnote{Supported in part by grants
R-155-000-035-112 and R-155-050-055-133/101 at National University of
Singapore} } \maketitle
%


\abstract{Suppose that $f(z)$ is a transcendental entire function
and that the Fatou set $F(f)\neq\emptyset$. Set
$$B_1(f):=\sup_{U}\frac{\sup_{z\in U}\log(|z|+3)}{\inf_{w\in
U}\log(|w|+3)}$$ and
$$B_2(f):=\sup_{U}\frac{\sup_{z\in U}\log\log(|z|+30)}{\inf_{w\in
U}\log(|w|+3)},$$ where the supremum $\sup_{U}$ is taken over all
components of $F(f)$. If $B_1(f)<\infty$ or $B_2(f)<\infty$, then
we say $F(f)$ is strongly uniformly bounded or uniformly bounded
respectively. In this article, we will show that, under some
conditions, $F(f)$ is (strongly) uniformly bounded.}

\vspace{0.3cm} \noindent {\bf Keywords: }Bounded Fatou components,
(multiply) connected component, (strongly) uniformly bounded

\vspace{0.3cm} \noindent {\bf AMS(2000) Subject classifications: }
Primary 37F10; Secondary 37F45

\section{Introduction}

Let $f(z)$ be a transcendental entire function, and denote by $f^n$, $n\in
\mathbb{N}$, the $n$th iterate of $f$. The set of normality, which is
the so called Fatou set $F(f)$, is defined to be the set of points, $z\in\mathbb{C}$,
such that the sequence $\{f^n\}$ is normal in some neighborhood of $z$.
The Fatou set $F(f)$ has good topological and analytical properties. For example,
$F(f)$ is open, and may contain $0$, $1$ or $\infty$ Fatou components.
If $F(f)$ contains only one Fatou component $D$, then $D$ is invariant
under $f$ and $D$ is unbounded. Another well-known set related to $f$ is
the Julia set, $J=J(f)=\mathbb{C}-F(f)$. $J(f)$ is a nonempty perfect set
which coincides with
$\mathbb{C}$, or is nowhere dense in $\mathbb{C}$.
 For the basic results in the dynamical theory of transcendental
functions, we refer the reader to the books \cite{hua1,mor1} and
the survey paper \cite{ber1}, and the papers of Fatou \cite{fat2}
and Julia \cite{jul1} for more about the iteration theory of
transcendental entire functions.

Throughout the paper, we shall use the following standard notation:
\begin{eqnarray*}
M(r,f)&=&\max\{|f(z)|:|z|=r\},\\
m(r,f)&=&\min\{|f(z)|:|z|=r\},\\
\lambda=\lambda(f)&=&\limsup_{r\rightarrow\infty}\frac{\log\log
M(r,f)}{\log r},\\
\rho=\rho(f)&=&\liminf_{r\rightarrow\infty}\frac{\log\log
M(r,f)}{\log r}.
\end{eqnarray*}
We call these maximum modulus, minimum modulus, order of $f$ and
lower order of $f$, respectively. If $E\subset [1,\infty)$, the
upper logarithmic density of the set $E$ is defined by
$$\overline{\log
dens}E=\limsup_{R\rightarrow\infty}\frac{1}{\log
R}\int_{E\cap(1,R)}\frac{dt}{t}.$$

\vspace{.2cm} \noindent Let $f(z)$ be a transcendental entire
function, then $f(z)$ has the following form
$f(z)=\sum_{j=0}^{\infty}a_jz^j$. We say that $f(z)$ has
\textbf{Fabry gaps} if
$$f(z)=\sum_{k=0}^{\infty}a_kz^{j_k} \qquad  \textrm{and} \qquad
\frac{j_k}{k}\rightarrow\infty$$ as $k\rightarrow\infty$.

\bd Let $f(z)$ be a transcendental entire function. If $f$
satisfies the following condition: \emph{there exist two positive
constants $\varepsilon_1$ and $\varepsilon_2$ with
$0<\varepsilon_1,\varepsilon_2<1$ such that \be\overline{\log
dens}E(f)\leq \varepsilon_1\label{001}\ee where \be
E(f)=\{r>1:\log m(r,f)\leq \varepsilon_2\log
M(r,f)\},\label{002}\ee}

\noindent then \textbf{we call $f\in\Delta$}. \ed

$\Delta$ contains many classes of transcendental entire functions
of any order; see examples 1, 2, 3 and 4 of the next section.

\vspace{.2cm}In the complex dynamical system theory, it is
interesting to investigate the boudedness of all Fatou components
of a function $f$. There is a vast extensive literature on the
boundedness of Fatou components, for example see
\cite{and1,bak1,hua1,rem1,sta1,wan1,yue1,yue2,zhe5} for the Fatou
sets of transcendental entire functions, \cite{zhe2,zhe3,zhe4,zhe5}
for the case of transcendental meromorphic functions.
In the above mentioned papers, the authors only considered whether or not
all  Fatou components of an entire or meromorphic transcendental function
are bounded. We could say that those results are coarse in some sense, and
can't describe the properties of Julia sets and Fatou sets perfectly. In
this paper we will show that all Fatou components of a class of
transcendental entire functions are not only bounded, but also
proportioned normally, that is, all  Fatou components have almost the
same ``size". As far as we know this is the first time to develop this property for
transcendental entire functions.

Let $f(z)$ be a transcendental entire function. Throughout the rest
of the paper, we assume that $F(f)\not=\emptyset$, and then $F(f)$
consists of one or countably many Fatou components. If $F(f)$
contains only one component $U$, then $U$ must be unbounded and
completely invariant under $f$. Furthermore,
 we need the following definition.

\bd Suppose that $f(z)$ is a transcendental entire function. Set
$$B_0(f):=\sup_{U}\frac{\sup_{z\in U}(|z|+1)}{\inf_{w\in
U}(|w|+1)},$$
$$B_1(f):=\sup_{U}\frac{\sup_{z\in U}\log(|z|+3)}{\inf_{w\in
U}\log(|w|+3)}$$ and
$$B_2(f):=\sup_{U}\frac{\sup_{z\in U}\log\log(|z|+30)}{\inf_{w\in
U}\log(|w|+3)},$$ where the supremum $\sup_{U}$ is taken over all
components of $F(f)$.

If $B_1(f)<\infty$ or $B_2(f)<\infty$, then we say the Fatou set
$F(f)$ is strongly uniformly bounded or uniformly bounded
respectively.\ed

\noindent {\bf Remark 1: }In the above definition, the numbers $3$
and $30$ are chosen in order to make $\log\log$ and $\log$ larger than 1,
They can be replaced by any other numbers greater than $e$ and
$e^e$ respectively. For any transcendental entire function $f(z)$,
if $F(f)=\emptyset$, then we put $B_1(f)=B_2(f)=1$, otherwise
$1<B_1(f),\ B_2(f)\leq\infty$. Also it is easy to see the
following facts:
\begin{enumerate}
\item $B_1(f)=B_2(f)=\infty$ if $F(f)$ contains unbounded
components; \item all components of $F(f)$ are bounded if $F(f)$ is
(strongly) uniformly bounded.
\end{enumerate}

\vspace{.2cm}\noindent {\bf Remark 2: }Strongly uniform boundedness
implies uniform boundedness. But the converse is not true; see example
4 of section 2 below.


\vspace{.4cm}The paper is organized as follows. In Section 2 we
present some examples of transcendental entire functions of any
order that belong to $\Delta$. In Section 3, we prove that the Fatou
sets are uniformly bounded for functions whose Fatou sets $F(f)$
contain at least two different obits of multiply connected
components or $\lambda(f^2)<\infty$. In Section 4, we show that the
same result holds for functions $f\in\Delta$ satisfying
\begin{enumerate}
\item there exist positive constants $M_0$, $C_1>1$, $C_2>1$ such that
$$[M_0,\infty)\subset \{r>0:\log M(C_1r,f)\geq C_2\log
M(r,f)\}.$$
\item or there exist positive constants $M_1$, $D_1>1$ and $D_2>1$ with $D_1>D_2^2$
such that
$$[M_1,\infty)\subset \{r>0:\log M(D_1r,f)< D_2\log M(r,f)\}.$$
\item or $0<\rho(f)\leq\lambda(f)<\infty$.
\end{enumerate}

\section{Examples}

\vspace{.2cm} \noindent  \textbf{Example 1: }(\cite{fuc1}) Let $f$
be an entire function of finite order with Fabry gaps. Then for any
$0<\varepsilon<1$, we have
$$\log m(r,f)>(1-\varepsilon)\log M(r,f)$$
holds for all $r$ outside a set of logarithmic density zero.
Therefore, the functions of this class belong to $\Delta$.

\vspace{.2cm}\noindent \textbf{Example 2: }( \cite{bar1}) Let $f$ be
a transcendental entire function of order $\lambda(f)<1/2$, and
suppose that $\lambda(f)<a<1/2$. Then if
$$E=\{r>0:\log m(r,f)\leq \cos(\pi a)\log M(r,f)\}$$
we have $\overline{\log dens} E\leq \lambda(f)/a,$ and so
$f\in\Delta$.

\vspace{.2cm}\noindent \textbf{Example 3: }Let $f$ be a
transcendental entire function of order $\lambda(f)<a<1/2$,
$n(>0)\in\mathbb{N}$ be an integer, and let $f_n(z)=f(z^n)$. If
$$E_n=\{r>0:\log m(r,f_n)\leq \cos(\pi a)\log M(r,f_n)\}.$$
Then $\overline{\log dens} E_n\leq \lambda(f)/a,$ and so
$f\in\Delta$.

\vspace{.2cm}\textbf{Note that: }$\lambda(f_n)=n\lambda(f)$ and
$\rho(f_n)=n\rho(f)$.

\vspace{.2cm}\noindent \textbf{Proof:} It is easy to note that, for
all $r>0$ \be m(r,f_n)=m(r^n,f)\textrm{ and
}M(r,f_n)=M(r^n,f).\label{ex01}\ee Set
$$E=\{r>0:\log m(r,f)\leq \cos(\pi a)\log M(r,f)\}$$
Then
$$r^n\in E \textrm{ if and only if }r\in E_n.$$

Now by combining the above Example 2 and (\ref{ex01}), we have
\begin{eqnarray*}
\overline{\log dens}E_n&=&\limsup_{r\rightarrow\infty}\frac{1}{\log
r}\int_{E_n\cap [1,r]}\frac{ds}{s}\ \ \ (\textrm{let }t=s^n) \\
&=&\limsup_{r\rightarrow\infty}\frac{1}{n\log
r}\int_{E\cap [1,r^n]}\frac{dt}{t}\\
&=&\overline{\log dens}E\\
&\leq&\frac{\lambda(f)}{a}.
\end{eqnarray*} \qed

\vspace{.2cm} \noindent {\bf Example 4: }Let $k_n(n\in\mathbb{N})$
denote any increasing sequence of positive integers, $r_1>2$ and
$C>0$ such that $0<C<(4e^2)^{-1}$. Suppose that $n_0$ is a
positive integer so that $2^{n_0-1}C>2r_1^{k_1}$. The sequence of
positive numbers $\{r_n\}_{n=1}^{\infty}$ is given inductively by
the equation: \be r_{n+1}
=C\left(1+\left(\frac{r_n}{r_1}\right)^{k_1}\right)\left(1+\left(\frac{r_n}{r_2}\right)^{k_2}\right)\cdots
\left(1+\left(\frac{r_n}{r_n}\right)^{k_n}\right).\label{003}\ee
Set
$$f(z)=C\prod_{n=1}^{\infty}\left(1+\left(\frac{z}{r_n}\right)^{k_n}\right).$$
This function has the following properties (see \cite{bak3} and
\cite{hua1}):
\begin{enumerate}
\item $F(f)$ contains multiply connected wandering domains and
thus every component of $F(f)$ is bounded; \item for sufficiently
large $n$, set $A_n=\{z:r^2_n\leq |z|\leq r_{n+1}^{1/2}\}$. Then
$f(A_n)\subset A_{n+1}$; \item for a given nonnegative number
$\lambda (0\leq\lambda\leq\infty)$, $\lambda(f)=\lambda$ if we take
$k_n=[r_n^{\lambda}]$, where $[a]$ denotes the largest integer that
is less than or equal to $a$. It can be verified that
$$\log M(r,f)=O(r^{\lambda})\textrm{ and }\log f(2r_n)>[r_n^{\lambda}]\log2.$$
\end{enumerate}

\vspace{.2cm}\textbf{Claim 1: }If $0<\lambda<\infty$, then $F(f)$ is
uniformly bounded.

\vspace{.2cm}\textbf{Proof of Claim 1: }Let $U_n$ be the multiply
connected Fatou components of $F(f)$ such that $A_n\subset U_n$,
$Q_1=\{z:|z|<r_1^{1/2}\}$, and
$Q_n=\{z:r_n^{1/2}<|z|<r_{n+1}^{2}\}$ for $n\geq 2$. It is easy to
check that there exists $n_0\in\mathbb{N}$ so that $r_{n+1}>2r_n$
for $n\geq n_0$. Without loss of generality, we may assume that
$r_{n+1}>2r_n$ for all $n\geq 1$. We note that the number of zeros
of $f(z)$ in $|z|\leq t$ $(r_n\leq t<r_{n+1})$ is
$n(t):=k_1+k_2+\cdot\cdot\cdot+k_n$, and there exists a positive
constant $C(\lambda)$ so that $n(t)\leq C(\lambda)r_{n}^{\lambda}$
for sufficiently large $n$ (see \cite{hua1}). Also we see from
(\ref{003}) that
$$r_{n+1}\leq Cr_{n}^{k_1+k_2+\cdot\cdot\cdot+k_{n-1}}\textrm{ and }r_{n}\leq
Cr_{n-1}^{k_1+k_2+\cdot\cdot\cdot+k_{n-2}}.$$ Then
$$r_{n+1}\leq C^{k_1+k_2+\cdot\cdot\cdot+k_{n-1}+1}r_{n-1}^{(k_1+k_2+\cdot\cdot\cdot+k_{n-1})^2}.$$
Thus
$$\frac{\log\log r_{n+1}}{\log r_{n-1}}<2\lambda+1$$
for sufficiently large $n$.

Now let $U$ be any Fatou component of $F(f)$. Then two cases should
be considered.

Case 1: $U=U_n$ for some $n\geq 1$.

In this case, it is easy to check that
$$U=U_n\subset\{z\in\mathbb{C}:r^{1/2}_{n-1}\leq|z|\leq r^{2}_{n+1}\}.$$
Thus
$$ \frac{\sup_{z\in U_n}\log\log(|z|+30)}{\inf_{w\in
U_n}\log(|w|+3)}\leq \frac{\log\log(r^2_{n+1}+30)}{\log(
r_{n-1}^{1/2}+3)} \leq\frac{6\log\log r_{n+1}}{\log r_{n-1}} \leq
6(2\lambda+1).$$

Case 2: $U\subset Q_n$ for some $n\geq 1$.

In this case, from the above argument, we also get that
$$ \frac{\sup_{z\in U}\log\log(|z|+30)}{\inf_{w\in
U}\log(|w|+3)}\leq 6(2\lambda+1).$$ Therefore, Claim 1 is proved.

\vspace{.2cm}\textbf{Claim 2: }$F(f)$ is uniformly bounded but not
strongly uniformly bounded.

\vspace{.2cm}\textbf{Proof of Claim 2: }It can be shown that for any
$m>1$, there exists a positive integer $n_m$ such that, for $n>n_m$,
$r_{n+1}>r_n^m.$ Then
$$\frac{\sup_{z\in U_n}\log(|z|+3)}{\inf_{w\in
U_n}\log(|w|+3)}\geq\frac{\log (r^{1/2}_{n+1}+3)}{\log
(r^2_n+3)}\geq \frac{1}{5}\frac{\log r_{n+1}}{\log
r_n}\geq\frac{m}{5}.$$ Therefore $$\sup_{U}\frac{\sup_{z\in
U}\log(|z|+3)}{\inf_{w\in U}\log(|w|+3)}=\infty.$$ Claim 2 also is
proved.

\vspace{.2cm}\textbf{Claim 3: }If all $k_n$ are odd then
$f(z)\in\Delta$.

\vspace{.2cm}\textbf{Proof of Claim 3: }For any given $r>0$, it is
clear that $m(r,f)=|f(-r)|$ and $M(r,f)=f(r)$. Set
$E_n=(4r_n,r_{n+1}/4)$ and $G_n=(r_n/4,4r_n)$ for sufficiently
large $n\geq 1$. Using a similar method as in the book
\cite{hua1}, we can also prove that there exists an integer
$N\in\mathbb{N}$ such that
$$M(r,f)=f(r)<4|f(-r)|=4m(r,f)$$
for any $n\geq N$ and for any $r\in E_n$. Choose a given number
$0<\varepsilon<1$. Set
$$G(f)=\{r>0:\log m(r,f)<(1-\varepsilon)\log M(r,f)\}.$$
Then, from the above arguments, we have
$$G(f)\subset\bigcup_{n=N}^{\infty}G_n.$$

Finally we need to show $\overline{\log dens}G(f)<1$. It is easy to see
that, for $n$ sufficiently large, $4r_n>16^{2n}.$ Hence
\begin{eqnarray*}
\overline{\log dens}G(f)&=&\limsup_{r\rightarrow\infty}\frac{1}{\log
r}\int_{G(f)\cap(1,r)}\frac{dt}{t}\\
&\leq&\limsup_{n\rightarrow\infty}\frac{1}{\log
4r_n}\int_{G(f)\cap(1,4r_n)}\frac{dt}{t}\\
&\leq&
\limsup_{n\rightarrow\infty}\frac{1}{\log(4r_{n})}\sum_{i=1}^{n}\int_{r_i/4}^{4r_{i}}\frac{dt}{t}\\
&=&\limsup_{n\rightarrow\infty}\left(\frac{n\log16}{\log 4r_{n}}\right)\\
&\leq&1/2.
\end{eqnarray*}
The above second inequality holds because  $G(f)$ is open
and $G(f)\bigcap(4r_n,r_{n+1}/4)=\emptyset$. \qed

\section{The Case of $f$ with Multiply Connected Components}

First we recall some definitions. A set $\cal X$ of Fatou
components of $F(f)$ is called an \textbf{orbit} if for any two
components $U_1,\ U_2\in\cal X$, there exist two positive integers
$m$ and $n$ so that $f^m(U_1)=f^n(U_2)$. For each domain $D$, we
will use \textbf{$\#(D)$} to denote the connectivity of $D$. Next,
we list some dynamical properties of transcendental entire functions
with multiply connected components, and others can be found in the
above mentioned books. It was I.N. Baker who proved the existence of
multiply connected wandering domains for transcendental entire
functions, see \cite{bak4,bak2,bak5,bak3}.

\vspace{0.2cm} \noindent {\bf Theorem A: }(\cite{bak1} or
\cite{hua1}, Theorem 4.1, Page 58) {\it Let $f$ be a
transcendental entire function and $D$ a multiply connected
component of $F(f)$. Then
\begin{enumerate}
\item All the Fatou components of $f$ are bounded.
\item $f^{m}(z)\to\infty$ uniformly on compact subsets of $D$.
\item For every Jordan curve $\gamma$ which is non-contractible in
$D$, index $(f^n(\gamma))\neq 0$ for all sufficiently large $n$.
\item $D$ must be wandering.
\end{enumerate}}

\vspace{0.2cm}Recently, in \cite{kis1}, Kisaka and Shishikura
proved the following results.

\vspace{0.2cm} \noindent {\bf Theorem B: }{\it
\begin{enumerate}
\item If $D$ is a multiply connected component of
$\infty$-connectivity, then $\#f^n(D)=\infty$ for all $n\geq1$.
\item If $\infty>\#D\geq2$, then there exists an integer $n_0\in\mathbb{N}$
such that $\#f^m(D)=2$ for all $m\geq n_0$.
\item For every $n\in\mathbb{N}$ with $n\geq2$ there exist a
transcendental entire function $f$ with a wandering domain $D$
with $\#D=n$ and $\#f^m(D)=2$ for every $m\geq1$.
\item There exists a
transcendental entire function with infinitely many grand orbits
of doubly connected wandering domains. That is, there exist doubly
connected wandering domains $D_i$ $(i\in\mathbb{N})$ such that if
$i\not= j$, then $f^m(D_i)\cap f^n(D_j)=\emptyset$ for any $m$,
$n\in\mathbb{N}$. Furthermore, this $f$ can be constructed so that
the Lebesgue measure of the Julia set $J(f)$ is positive.
\end{enumerate}}

Now we state our first result.

\vspace{.2cm}\noindent {\bf Theorem 1: }{\it Let $f(z)$ be a
transcendental entire function with $F(f)$ containing multiply
connected components.
\begin{enumerate}
\item If $\lambda(f^2)<\infty$, then $F(f)$ is uniformly bounded.
\item Suppose that $F(f)$ contains at least two different obits $\cal
X$ and $\cal Y$ which consist of only multiply connected
components, suppose further that $\lambda(f)<\infty$. Then $F(f)$
is uniformly bounded.
\item $B_0(f)=\infty.$
\end{enumerate}}

\subsection{Proof of the Theorem 1}

Since $F(f)$ has multiply connected components, we may assume that
$U$ is a multiply connected component of $F(f)$. Let $U_n(n\geq
0)$ be the Fatou components so that $f^n(U)\subset U_n$ and $U_0=U$.
There exists a simple closed curve $\gamma$ in $U$ which is not
null-homotopic in $U$. Define $\gamma_{n}=f^n(\gamma)$. Since $U$ is
multiply connected, $f^n|_U\rightarrow\infty$ as
$n\rightarrow\infty$. By a theorem of Baker (\cite{bak1} or
\cite{hua1}), $0\in int(\gamma_n)$ for sufficiently large $n$. And
by Lemma 7 of Bergweiler \cite{ber1}, there exist a constant $K>1$
and a positive integer $n_0$ such that $|f^n(z_1)|\leq |f^n(z_2)|^K$
for all $z_1$, $z_2\in\gamma$ and $n\geq n_0$. Let
$r_n=\min\{|f^n(z)|:z\in\gamma\}$. Then $\gamma_n\subset
ann(r_n,r_n^K)$ for $n\geq n_0$, where $ann(r,R)$ denotes the
annulus around $0$ with radii $r$ and $R$.

Note also that the sequence $\{r_n\}$ satisfies
$\lim_{n\rightarrow\infty}r_n=\infty$ and \be
\min_{z\in\gamma_n}|f(z)|=
r_{n+1}\geq(\max_{z\in\gamma_n}|f(z)|)^{1/K}\geq
M(r_n,f)^{1/K}\label{401}\ee for $n\geq n_0$.


\vspace{.2cm}\noindent \textbf{Case 1:} If $\lambda(f^2)<\infty$, we
show that $F(f)$ is uniformly bounded.

\vspace{.2cm}It is easy to see that \be
U_n\subset\{z\in\mathbb{C}:r_{n-1}\leq|z|\leq
r^K_{n+1}\}\label{402}\ee and that we can assume $U_n$ is doubly
connected for $n\geq n_0$ by Theorem B. Now let $\Gamma_n$ be the
outer boundary of $U_n$ for $n\geq n_0$. Then we know that in order
to prove this case, we only need to consider the components $E$ that
lie outside of $\Gamma_{n_0}$, and such components can be separated
into the following two classes:

\textbf{Class 1: }$E=U_n$ for some $n>n_0$;

\textbf{Class 2: }$E$ lies between the outer boundaries of $U_n$ and
$U_{n+1}$ for some $n\geq n_0$;


\textbf{Claim 1: }There exists a positive number $D_5$ so that
$$\sup_{n\geq n_0}\frac{\sup_{z\in U_n}\log\log(|z|+30)}{\inf_{w\in
U_n}\log(|w|+3)}<D_5.$$

\textbf{Proof of Claim 1: }It's easy to check that, for $n\geq n_0$,
$$r_{n+1}=\min_{z\in\gamma}|f^{n+1}(z)|=\min_{z\in\gamma_{n-1}}|f^2(z)|\leq
M(r_{n-1}^K,f^2).$$ Combining this, $\lambda(f^2)<\infty$ and
(\ref{402}), we get
\begin{eqnarray*}
\frac{\sup_{z\in U_n}\log\log(|z|+30)}{\inf_{w\in U_n}\log(|w|+3)}
&\leq&\frac{\log\log(r^K_{n+1}+30)}{\log(r_{n-1}+3)} \\
&\leq& K\frac{\log\log r_{n+1}}{\log r_{n-1}}\\
&\leq& K^2\frac{\log\log M(r_{n-1}^K,f^2)}{\log r^K_{n-1}}\\
&\leq& K^2(\lambda(f^2)+1)\\
&<&\infty
\end{eqnarray*}
for sufficiently large $n$. Now Claim 1 follows.

\textbf{Claim 2: }There exists a positive number $D_6$ so that
$$\sup_{E}\frac{\sup_{z\in E}\log\log(|z|+30)}{\inf_{w\in
E}\log(|w|+3)}\leq D_6,$$ where $\sup_E$ is taken over all
components in Class 2.

\textbf{Proof of Claim 2: }Let $E$ be such a component so that $E$
lies between the outer boundaries of $U_n$ and $U_{n+1}$ for some
$n\geq n_0$. Then it is easy to check that
$$r_n\leq|z|\leq r^K_{n+1}\textrm{ for }z\in E$$
and that \be
r_{n+1}=\min_{z\in\gamma}|f^{n+1}(z)|=\min_{z\in\gamma_{n}}|f(z)|\leq
M(r_{n}^K,f).\label{444}\ee Hence
\begin{eqnarray*}
\frac{\sup_{z\in E}\log\log(|z|+30)}{\inf_{w\in E}\log(|w|+3)}
&\leq& \frac{\log\log(r^K_{n+1}+30)}{\log(r_{n}+3)}\\
&\leq& K\frac{\log\log r_{n+1}}{\log r_{n}}\\
&\leq& K^2\frac{\log\log M(r_{n}^K,f)}{\log r^K_{n}}\\
&\leq& K^2(\lambda(f)+1)
\end{eqnarray*}
for sufficiently large $n$. Now Claim 2 and furthermore Case 1
follow.

%

\vspace{.2cm}\noindent \textbf{Case 2:} Suppose that $\cal X$ and
$\cal Y$ are two different orbits of multiply connected Fatou
components of $F(f)$. We show that $F(f)$ is uniformly bounded.

\vspace{.2cm}Without loss of generality, we may assume that
$U_n(n\geq 0)\in\cal X$ and that $U'_n(n\geq 0)\in\cal Y$ is another
sequence of multiply connected components with $f^n(U'_0)\subset
U'_n$. We also assume that $\gamma'$, $\gamma'_n=f^n(\gamma')(n\geq
1)$ and $r'_n(n\geq 1)$ are defined in the same way as $\gamma$,
$\gamma_n=f^n(\gamma)(n\geq 1)$ and $r_n(n\geq 1)$, respectively.

Now, fix a sufficiently large $n_0\in\mathbb{N}$. Then all the
Fatou components of $F(f)$ can be put in the following three
classes:
\begin{eqnarray*}
\cal A&=&\{U\subset F(f):U\subset \textrm{int}(\gamma_{n_0})\};\\
\cal B&=&\{U\subset F(f):U=U_n\textrm{ for some }n\geq n_0\};\\
\cal C&=&\{U\subset F(f):U\subset \textrm{int}(\gamma_n)\textrm{ and
}U\subset \textrm{out}(\gamma_{n-1})\textrm{ for some }n>n_0\}.
\end{eqnarray*}

\textbf{Subcase 1:} It is trivial to see that there exists a
positive constant $C_3$ so that
$$\frac{\sup_{z\in U}\log\log (|z|+30)}{\inf_{w\in
U}\log(|w|+3)}\leq C_3$$ for any $U\in\cal A$.

\textbf{Subcase 2:} For any given $U_{n}\in \cal B$, it follows from
the fact, $\cal X$ and $\cal Y$ are different great orbits, that
there exists some $m\in\mathbb{N}$ so that
$$U_{n}\subset (r'_{m},r'^K_{m+1}).$$
Clearly, $m\rightarrow\infty$ as $n\rightarrow\infty$. Since
$$
r'_{m+1}=\min_{z\in\gamma'}|f^{m+1}(z)|=\min_{z\in\gamma'_m}|f(z)|\leq
M(r'_m,f),$$ we have for large $n$
\begin{eqnarray*}
\frac{\sup_{z\in_{U_{n}}}\log\log(|z|+30)}{\inf_{w\in
U_{n}}\log(|w|+3)}&\leq&
\frac{\log\log(r'^K_{m+1}+30)}{\log(r'_m+3)}\\
&\leq&K^2\cdot\frac{\log\log M(r'^K_m,f)}{\log r'^K_m}\\
&\leq& K^2(\lambda(f)+1)\\
&<&\infty.
\end{eqnarray*}

\textbf{Subcase 3:} For any $U_{n}\in\cal C$, similarly, for large
$n$, we have by (\ref{444}) that
\begin{eqnarray*}
\frac{\sup_{z\in_{U_{n}}}\log\log(|z|+30)}{\inf_{w\in
U_{n}}\log(|w|+3)}&\leq&
\frac{\log\log(r^K_{n}+30)}{\log(r_{n-1}+3)}\\
&\leq&K^2\cdot\frac{\log\log M(r^K_{n-1},f)}{\log r^K_{n-1}}\\
&\leq& K^2(\lambda(f)+1)\\
&<&\infty.
\end{eqnarray*}
Case 2 follows.

\vspace{.2cm}\noindent \textbf{Case 3:} Let $U$ be a multiply
connected Fatou component of $F(f)$. In \cite{zhe6}, Zheng proved
that if $n\in\mathbb{N}$ is sufficiently large then there exists a
round annulus $D_n=\{r_n <|z|<R_n\}\subset f^n(U)$ such that
$dist(0;D_n)\rightarrow\infty$ and $\mod(D_n)\rightarrow\infty$ as
$n\rightarrow\infty$. Therefore
$$\sup_{D}\frac{\sup_{z\in D}(|z|+1)}{\inf_{w\in
D}(|w|+1)}=\infty,$$ where the supremum is taken over all Fatou
components of $F(f)$. Thus Case 3 follows and therefore the proof of
Theorem 1 is complete.\qed

\section{The Case of $f\in\Delta$}

\vspace{.2cm}\noindent {\bf Theorem 2: }{\it Suppose that
$f(z)\in\Delta$ is a transcendental entire function and that there
exist positive constants $M_0$, $C_1>1$, $C_2>1$ such that \be
[M_0,\infty)\subset \{r>0:\log M(C_1r,f)\geq C_2\log
M(r,f)\},\label{004}\ee then $F(f)$ is strongly uniformly bounded.}

\vspace{.2cm}\noindent {\bf Theorem 3: }{\it Suppose that $f(z)$ is
a transcendental entire function with $\lambda(f)<1/2$. If there
exist positive constants $M_1$, $D_1>1$ and $D_2>1$ with $D_1>D_2^2$
such that
$$[M_1,\infty)\subset \{r>0:\log M(D_1r,f)< D_2\log M(r,f)\},$$
then $F(f)$ is uniformly bounded.}

\vspace{.2cm}\noindent {\bf Theorem 4: }{\it Suppose that
$f(z)\in\Delta$ is a transcendental entire function. If
$0<\rho(f)\leq\lambda(f)<\infty$, then $F(f)$ is uniformly
bounded.}

\vspace{.2cm}It's natural to post the following two conjectures.

\vspace{.2cm}\noindent {\bf Conjecture 1: }{\it If $f(z)\in\Delta$
and $F(f)\not=\emptyset$, then $F(f)$ is uniformly bounded.}

\vspace{.2cm}\noindent {\bf Conjecture 2: }{\it If $f(z)$ is a
transcendental entire function and $F(f)\not=\emptyset$, then
$B_0(f)=\infty$.}

\subsection{Proof of Theorem 2}

First we recall some facts:

\textbf{Fact 1 }(see \cite{zhe6}): Let $f(z)$ be a transcendental
entire function. If there exist two constants $C_1>1$ and $C_2>1$ so
that
$$
\log M(C_1r,f)\geq C_2\log M(r,f)$$ for sufficiently large $r$,
then $F(f)$ contains no multiply connected Fatou components.

\textbf{Fact 2 }(see \cite{car1}): Let $U$ be a simply connected
hyperbolic domain of $\mathbb{C}$, $\rho_U(z)$ the density of the
hyperbolic metric on $U$ and $\lambda(z_1,z_2)$ the hyperbolic
distance between $z_1$ and $z_2$ on $U$. Then
$$\frac{1}{2d(z,\partial U)}\leq\rho_U(z)\leq\frac{2}{d(z,\partial
U)},$$ where $d(z,\partial U)$ is the Euclidean distance between
$z\in U$ and $\partial U$.

\vspace{.2cm}\textbf{Claim 1: }There exist two positive numbers
$M_1(>M_0)$ and $h>1$ so that, for all $r>M_1$, there exists
$r'\in(r,r^h)$ satisfying \be m(r',f)>M(r,f)^h.\label{005}\ee

\vspace{.2cm}\textbf{Proof of Claim 1: }We recall that $f$
satisfies the following condition: there exist two positive
constants $\varepsilon_1$ and $\varepsilon_2$ with
$0<\varepsilon_1,\varepsilon_2<1$ such that \be\overline{\log
dens}E(f)\leq \varepsilon_1\label{3.01}\ee where \be
E(f)=\{r>1:\log m(r,f)\leq \varepsilon_2\log
M(r,f)\}.\label{3.02}\ee Now we choose a positive number $c$ with
$\varepsilon_1>c\varepsilon_1+\varepsilon_1^2$.

\vspace{.2cm} \textbf{Claim 1.1: }There exist two positive numbers
$D$ and $r_0$ such that, for any $r\geq r_0$, there must be some
number $s$ with
$$s\in(r,r^{(1+\varepsilon_1/c)})
\textrm{ and } m(s,f)\geq M(r,f)^{D}.$$

\textbf{Proof of Claim 1.1: }We claim that there exists a positive
number $r_0$ such that, if $R\geq r_0$, there is some value
$r_R\in (R,R^{1+\varepsilon_1/c})$ so that
$$\log m(r_R,f)>\varepsilon_2\log M(r_R,f).$$
Suppose this is not true. Then there exists a sequence $\{R_j\}$,
$R_j\rightarrow\infty$ as $j\rightarrow\infty$ such that
$$(R_j,R_j^{1+\varepsilon_1/c})\subset E(f).$$
Thus
$$\overline{\log
dens}E(f)\geq\limsup_{j\rightarrow\infty}\frac{1}{(1+\varepsilon_1/c)\log
R_j}\int_{R_j}^{R_j^{1+\varepsilon_1/c}}\frac{dt}{t}=\frac{\varepsilon_1}{\varepsilon_1+c}>\varepsilon_1.$$
This obviously contradicts (\ref{3.01}).

Then, for sufficiently large $r$, there exists $s$ with $s\in
(r,r^{(1+\varepsilon_1/c)})$ such that \be\log
m(s,f)>\varepsilon_2\log M(s,f).\label{3.03}\ee

Now applying Hadamard's three-circles theorem to the three circles
$|z|=1,r,s$, we get that

\be\log M(s,f)\geq \log M(r,f)-\frac{\varepsilon_1}{c}\log
M(1,f).\label{3.04}\ee Clearly, when $r$ is chosen sufficiently
large, we have

\be\frac{\varepsilon_1}{c}\log M(1,f)\leq\frac{1}{2}\log
M(r,f).\label{3.05}\ee By combining (\ref{3.03}), (\ref{3.04}) and
(\ref{3.05}), we see that there exists a positive constant $r_0$
such that, when $r>r_0$, we have
$$
\log m(s,f)\geq\frac{\varepsilon_2}{2}\log M(r,f).
$$
Therefore, we have proved Claim 1.1 with $D=\varepsilon_2/2$.

\vspace{.2cm} \textbf{Claim 1.2: }There exist two positive numbers
$M_1(>M_0)$ and $h>1$ so that, for all $r>M_1$, there exists
$r'\in(r,r^h)$ satisfying
$$ m(r',f)>M(r,f)^h.$$

\textbf{Proof of Claim 1.2: }From (\ref{004}) one easily gets that
\be \log M(r_2,f)\geq C_2^{\frac{\log r_2-\log r_1}{\log
C_1}-1}\log M(r_1,f)\label{204}\ee where $r_2>r_1>M_0$.

Now let $\beta=1+\varepsilon_1/c$. From Claim 1.1, there exists
$r'\in(r^{\beta^2},r^{\beta^3})$ such that
$$m(r',f)>M(r^{\beta^2},f)^{D}.$$
It is easy to check that there exists a positive number
$r'_0(>r_0)$ such that, for any $r\geq r'_0$,
$$h:=DC_2^{\frac{(\beta^2-\beta)\log r}{\log C_1}-1}>\beta^2(>1).$$
Then by (\ref{204}), if $r\geq r'_0$,
$$M(r^{\beta^2},f)>M(r^{\beta},f)^{h/D}.$$
The above two inequalities imply that
$$m(r',f)>M(r^{\beta},f)^h$$
for some $r'\in(r^{\beta},r^{\beta^3})$, where $r>r'_0$. Claim 1.2
follows, and therefore Claim 1 holds.

\vspace{.2cm}In order to prove the Theorem, we will show that
there exists a positive $D$ so that

$$\sup_U\frac{\sup_{z\in U}\log(|z|+3)}{\inf_{w\in
U}\log(|w|+3)}<D,$$ where the supremum is taken over all Fatou
components of $F(f)$. To do this, we need to prove the following
Claims 2 and 3.

\vspace{.2cm}\textbf{Claim 2: }There exists a positive number
$D_1$ so that
$$\sup_U\frac{\sup_{z\in U}\log(|z|+3)}{\inf_{w\in
U}\log(|w|+3)}<D_1,$$ where the supremum is taken over all simply
connected components of $F(f)$.

\vspace{.2cm}\textbf{Proof of Claim 2: }If Claim 2 does not hold,
then there exists a simply connected component $U$ and a point
$z_0\in U$ so that $|z_0|\geq M_1$ and \be\frac{\sup_{z\in
U}\log(|z|+3)}{\log(|z_0|+3)}>2h.\label{e00}\ee

Suppose that $U_n(n\geq 1)$ is a simply connected component so that
$f^n(U)\subset U_n$. Then $U_n(n\geq 1)$ is also simply connected.
Now by Claim 1 and the fact that $J(f)\not=\emptyset$, we can take
$z_0,z_1\in U$, $R'_0>0$, $R_0>M_1$ and a Jordan curve
$\gamma_0\subset U$ so that
\begin{enumerate}

\item $z_0,z_1\in \gamma_0$, $R'_0\in(R_0,R_0^h)$, $|z_0|=R_0$ and
$|z_1|=R'_0$; \item $m(R'_0,f)>M(R_0,f)^h$; \item
\be\{z:|z|<R_0\}\cap J(f)\not=\emptyset.\label{e01}\ee
\end{enumerate}
This is possible because of (\ref{e00}). It follows from (\ref{e01})
that \be M(R_0,f^n)\rightarrow\infty\textrm{ as
}n\rightarrow\infty.\label{e02}\ee

First, it is trivial to check that
\begin{enumerate}
\item $|f(z_1)|\geq m(R'_0,f)>M(R_0,f)^h$;\item
$|f(z_1)|\geq|f(z_0)|^h.$
\end{enumerate}

In the following we will mainly use Zheng-Wang's idea (see
\cite{zhe5}) to complete our argument.

Now put $R_1=M(R_0,f)$, $\gamma_n=f^n(\gamma_0)$ for $n\geq 1$,
$\Gamma_1=\{z:|z|=R_1\}$ and $\Gamma_1^h=\{z:|z|=R_1^h\}$. Then by
the above argument we know that $\gamma_1\cap\Gamma_1\not=\emptyset$
and $\gamma_1\cap\Gamma_1^h\not=\emptyset$. Thus by Claim 1 there
are a positive number $R'_1$ and a point $z_2\in\gamma_0$ so that
$$R_1<R'_1<R_1^h,\ f(z_2)\in\{z:|z|=R'_1\}$$
and that
\begin{eqnarray*}
|f^2(z_2)|&=&|f(f(z_2))|\geq m(R'_1,f)\geq M(R_1,f)^h\\
&=&M(M(R_0,f),f)^h\geq M(R_0,f^2)^h\\&\geq& |f^2(z_0)|^h.
\end{eqnarray*}
That is $|f^2(z_2)|\geq \max\{M(R_0,f^2)^h, |f^2(z_0)|^h\}.$
Similarly, by putting $R_2=M(R_1,f)$, $\Gamma_2=\{z:|z|=R_2\}$ and
$\Gamma_2^h=\{z:|z|=R_2^h\}$, we can choose a point $z_3\in
\gamma_0$ so that $$|f^3(z_3)|\geq
\max\{M(R_0,f^3)^h,|f^3(z_0)|^h\}.$$ Repeating the above process
inductively, we can deduce that there is a point $z_n\in\gamma_0$ so
that \be |f^n(z_n)|\geq \max\{M(R_0,f^n)^h,|f^n(z_0)|^h\}\textrm{
for }n\geq 1.\label{e03}\ee Now, for any $n$, pick points
$a_n\in\partial U_n$ so that $|a_n|\leq |f^n(z_0)|$, and we have by
Fact 2,
$$
\rho_{U_n}(z)\geq\frac{1}{2|z-a_n|}\geq\frac{1}{2(|z|+|a_n|)}.$$
Then \be \rho_{U_n}(f^n(z_0),f^n(z_n))\geq
\int^{|f^n(z_n)|}_{|f^n(z_0)|}\frac{dr}{2(r+|a_n|)}\geq\frac{1}{2}
\log\frac{|f^n(z_n)|+|a_n|}{|f^n(z_0)|+|a_n|}.\label{e05}\ee It is
obvious that \be
\rho_{U_n}(f^n(z_0),f^n(z_n))\leq\rho_U(z_0,z_n)\leq
\lambda(\gamma_0)<+\infty\label{e06}.\ee Here $\lambda(\gamma_0)$
is the hyperbolic length of $\gamma_0$ in $U$.

Combining (\ref{e03}), (\ref{e05}) and (\ref{e06}), we get that \be
2|f^n(z_0)|e^{2\lambda(\gamma_0)}\geq(|f^n(z_0)|+|a_n|)e^{2\lambda(\gamma_0)}\geq
|f^n(z_n)|+|a_n|\geq |f^n(z_0)|^h.\label{e07}\ee

Since $U$ is a simply connected component, by (\ref{e02}) and
(\ref{e03}) we have that $|f^n(z_n)|\rightarrow\infty$ as
$n\rightarrow\infty$. Thus $|f^n(z)|\rightarrow\infty$ as
$n\rightarrow\infty$ for any $z\in U$, and it means that
$|f^n(z_0)|\rightarrow\infty$ as $n\rightarrow\infty$. This
contradicts to (\ref{e07}) with $h>1$. Therefore Claim 2 follows.

\vspace{.2cm}\textbf{Claim 3: }There exist no multiply connected
components in the Fatou set $F(f)$.

\vspace{.2cm}\textbf{Proof of Claim 3: }This follows from Fact 1.

\vspace{.2cm}Theorem 2 now follows from Claims 2 and 3. \qed

\subsection{Proof of Theorem 3}

\bl \label{lem2.1}Suppose that $f(z)$ is a transcendental entire
function, and that there exist sequences $\{R_n\}$,
$\{t_n\}\rightarrow\infty$ and $\{c(n)\}$ with $c(n)>1$ $(n\geq
1)$ and $a=\sup\{c(n)\}<\infty$ so that:
\begin{enumerate}
\item $M(R_n,f)=R_{n+1}$; \item $R_n< t_n< R_{n}^{c(n)}$; \item
$m(t_n,f)>R_{n+1}^{c(n+1)}$ for all $n\geq 1$.
\end{enumerate}
Then there is no Fatou component $U$ so that $$U\cap
\{z:|z|=R_n\}\not=\emptyset\textrm{ and }U\cap
\{z:|z|=R_{n+1}^b\}\not=\emptyset$$ for some $n\geq 1$ and any
$b>a$. \el

\noindent {\bf Proof of Lemma \ref{lem2.1}: }Without loss of
generality, we may assume that $J(f)\cap
\{z:|z|<R_1\}\not=\emptyset$. Then $M(R_1,f^n)\rightarrow\infty$ as
$n\rightarrow\infty$. Thus $R_{n}\rightarrow\infty$ as
$n\rightarrow\infty$.

This Lemma is similar to a result of Baker, see \cite{hua1} or
\cite{bak1}. Set
$$C_n=\{z:|z|=R_n\},C_n^1=\{z:|z|=R_n^{c(n)}\} \textrm{ and }
C_n^2=\{z:|z|=t_n\}$$ for $n\geq 1$.

Assume that there exists a Fatou
component $U$ so that $U\cap C_1\not=\emptyset$ and $U\cap
C_2^1\not=\emptyset$. We need to deduce a contradiction.

It is trivial to see that $U$ contains a path $\Gamma$ joining
points $w_1\in C_1$ and $w_{2}^1\in C_{2}^1$, and obviously $\Gamma$
contains a point $w_2^2\in C_2^2$. Then $|f(w_1)|\leq R_{2}$ and
$|f(w_{2}^2)|> R_{3}^{c(3)}$ from the conditions of this Lemma.
Hence $f(\Gamma)$ must contain an arc joining a point $w_{2}\in
C_{2}$ to a point $w_{3}^1\in C_{3}^1$. This process can be
repeated, and $f^k(U)$ contains an arc of $f^k(\Gamma)$ which joins
a point $w_{k+1}\in C_{k+1}$ to a point $w_{k+2}^1\in C_{k+2}^1$.
Thus, on $\Gamma$, the function $f^k$ takes a value of modulus at
least $R_{k+1}\rightarrow\infty$ as $k\rightarrow\infty$. Then we
conclude that $f^k\rightarrow\infty$ as $k\rightarrow\infty$ locally
uniformly on $U$. It follows that there exists $n_1\in\mathbb{N}$
such that, for $k\geq n_1$, we have $|f^k(z)|>1$ for all $z\in
\Gamma$. Now by a result of \cite{ber1}, there exists a positive
constant $B$ such that $|f^k(w)|\leq |f^k(z)|^{B}$ for all $w,z\in
\Gamma$ and for all $k\geq n_1$. Now for any $k\geq n_1$, we pick
points $z_k,z_k^1\in \Gamma$ such that $f^k(z_k)=w_{k+1}$ and
$f^k(z_k^1)=w_{k+2}^1$. Thus, for each $k>n_1$, we have
$$M(R_{k+1},f)=R_{k+2}<R_{k+2}^{c(k+2)}=|f^k(z_k^1)|<|f^k(z_k)|^{B}=
R_{k+1}^{B}.$$ This contradicts the fact that $f$ is a
transcendental entire function. \qed

\vspace{.2cm}\noindent\textbf{Proof of Theorem 3: }Since $f(z)$
satisfies the conditions of the theorem, by a result of Hua and Yang
(see \cite{hua2}, Lemma 2), for any given sufficiently large $R_1$
($>\min_{z\in J(f)}|z|$) there exists a sequence $\{t_n\}$ so that,
for $n\geq 1$, \be M(R_n,f)=R_{n+1}\label{501};\ee \be
R_n<R_n^{2+\frac{2}{2n+1}}< t_n<
R_{n}^{2+\frac{1}{n}}\label{502};\ee \be\log
m(t_n,f)>\left(1-\frac{2}{(2n+1)^3}\right)\log
M(R_{n}^{2+\frac{2}{2n+1}},f)\label{503}.\ee

Now applying Hadamard's three-circles theorem to $|z|=1$, $R_n$ and
$R_n^{2+\frac{2}{2n+1}}$, we get
$$\log M(R_n£¬f)\leq \log
M(1,f)+\frac{1}{2+\frac{2}{2n+1}}\log
M(R_n^{2+\frac{2}{2n+1}},f).$$

This implies the following,
\begin{eqnarray}
& &\left[1-\frac{1}{(2n+1)^3}\right]\log
M(R_n^{2+\frac{2}{2n+1}},f)\nonumber\\
&>&\left[2+\frac{1}{n+1}\right]\log
M(R_n,f)\nonumber\\
&
&+\left[\left(2+\frac{2}{2n+1}\right)\left(1-\frac{1}{(2n+1)^3}\right)-
\left(2+\frac{1}{n+1}\right)\right]\log M(R_n,f)\nonumber\\
&
&-\left(2+\frac{2}{2n+1}\right)\left(1-\frac{1}{(2n+1)^3}\right)\log
M(1,f)\nonumber\\
&=&\left[2+\frac{1}{n+1}\right]\log
M(R_n,f)+\frac{8n^3+8n^2-2n-3}{(2n+1)^4(n+1)}\log
M(R_n,f)\nonumber\\
&
&-\left(2+\frac{2}{2n+1}\right)\left(1-\frac{1}{(2n+1)^3}\right)\log
M(1,f).\label{504}
\end{eqnarray}
Since $f$ is transcendental, we may suppose that, for all $n\geq
1$, \be\frac{8n^3+8n^2-2n-3}{(2n+1)^4(n+1)}\log
M(R_n,f)>\left(2+\frac{2}{2n+1}\right)\left(1-\frac{1}{(2n+1)^3}\right)\log
M(1,f).\label{505}\ee

It's easy to see that (\ref{503}), (\ref{504}) and (\ref{505})
give the following result, \be
m(t_n,f)>R_{n+1}^{2+\frac{1}{n+1}}.\label{506}\ee

We fix $R_1$, and set
$$\triangle_1=\{z\in\mathbb{C}:|z|\leq R_2^3\}
\textrm{ and }\triangle_2=\{z\in\mathbb{C}:|z|\geq R_1\}.$$ Then all
the Fatou components can be put into the following three classes
(there could be a component in both of the first two classes):
\begin{eqnarray*}
\cal A&=&\{U\subset F(f):U\subset \triangle_1\},\\
\cal B&=&\{U\subset F(f):U\subset \triangle_2\},\\
\cal C&=&\{U\subset F(f):U\cap
\partial(\triangle_1)\not=\emptyset\textrm{ and
}U\cap\partial(\triangle_2)\not=\emptyset\}.
\end{eqnarray*}
By Lemma \ref{lem2.1} and the above arguments, we easily see that
$\cal C=\emptyset$. To prove that $F(f)$ is uniformly bounded, we
only need to consider the Fatou component in $\cal B$.

Let $U'$ be any Fatou component of $F(f)$. If
$U'\cap\{w:|w|=M(|z|,f)\}=\emptyset$ for any $z\in U'$, then we
put $\Theta_2(U')=1$. Otherwise, set
$$\Theta_2(U')=\sup\{t\geq 1:U'\cap\{w:|w|=M(|z|,f)^t,\textrm{ where }z\in
U'\}\neq\emptyset\}.$$ Set
$$\Theta_2=\sup_{U'\in\cal B}\Theta_2(U').$$

In the following, we show that $\Theta_2\leq 3$. Suppose on the
contrary that there exists $U_0\in\cal B$ such that
$\Theta_2(U_0)>3$. Then we can choose a point $z_1\in U_0$ so that
$$U_0\cap\{z:|z|=T_2^3\}\not=\emptyset$$
where $T_2=M(T_1,f)$ and $T_1=|z_1|$.

By (\ref{501}), (\ref{502}) and (\ref{506}), from this $T_1$, we
can get two sequences $\{T_n\}$ and $\{t_n\}$ such that
\begin{enumerate}
\item $M(T_n,f)=T_{n+1}$; \item $T_n< t_n< T_{n}^{2+\frac{1}{n}}$;
\item $m(t_n,f)>T_{n+1}^{2+\frac{1}{n+1}}$
\end{enumerate}
for $n\geq 1$. This obviously contradicts Lemma \ref{lem2.1}. So
$\Theta_2(U_0)\leq 3$. Thus
$$\sup_{U'\in\cal B}\frac{\sup_{w\in U'}\log\log(|w|+30)}{\inf_{z\in U'}\log(|z|+3)}\leq
2\lambda(f)<\infty.$$ It's easy to see that there exists a number
$M>0$ so that
$$\sup_{U\subset F(f)}\frac{\sup_{w\in U}\log\log(|w|+30)}{\inf_{z\in U}\log(|z|+3)}\leq
M.$$ Theorem 3 follows. \qed

\subsection{Proof of Theorem 4}

\vspace{.2cm} \textbf{Claim 1: }There exist two positive numbers
$D$ and $r_0$ such that, for any $r\geq r_0$, there must be some
number $s$ with
$$s\in(r,r^{D})
\textrm{ and } m(s,f)\geq M(r,f)^{D}.$$

\textbf{Proof of Claim 1: }This follows from the Claim 1 of Theorem
2.

\vspace{.2cm} \textbf{Claim 2: }Now take two positive constants $A$
and $B$ with $A>\lambda(f)/\rho(f)$ and $B>1$. For any
sufficiently large $R_1$, there exist a sequence $\{R_n\}$,
$R_n\rightarrow\infty$ as $n\rightarrow\infty$ such that \be
M(R_n^A,f)>M(R_n,f)^B\label{3.06}\ee for sufficiently large $n$.

\textbf{Proof of Claim 2: }Choose any $R>\min_{z\in J(f)}|z|$. Since
$J(f)\not=\emptyset$ for any transcendental entire function $f$, we
have $0<R<\infty$. Since $J(f)\cap \{z:|z|<R\}\not=\emptyset$,
$f^n(\{z:|z|<R\})$ must contain any given compact set of
$\mathbb{C}$ for sufficiently large $n$ as long as the compact set
doesn't meet an exceptional point. So if we set $R_1=M(R,f)$,
$R_{n+1}=M(R_n,f)(n=1,2,3,\cdots.)$, then $R_n\rightarrow\infty$ as
$n\rightarrow\infty$. We claim that this sequence is what we want.

Suppose on the contrary that (\ref{3.06}) does not hold. Then
there must be a subsequence $\{R_{n_j}\}\subset \{R_n\}$ such that
$$M(R_{n_j}^A,f)\leq M(R_{n_j},f)^B$$
for all $j\geq 1$. Then
$$A\cdot\frac{\log\log M(R^{A}_{n_j},f)}{\log R^{A}_{n_j}}\leq
\frac{\log B+\log\log M(R_{n_j},f)}{\log R_{n_j}}.$$ Thus
$$A\rho(f)\leq \lambda(f),$$
which contradicts the choice of $A$.

\vspace{.2cm}\textbf{Claim 3: }$F(f)$ is uniformly bounded.

\textbf{Proof of Claim 3: }Take any sufficiently large
$R(>\min_{z\in J(f)}|z|$). Now choose a sequence $\{R_n\}$: \be
R_1=M(R,f),\ R_{n+1}=M(R_n,f),\ n=1,2,3,\cdots,\
R_n\rightarrow\infty\textrm{ as }n\rightarrow\infty.\label{3.07}\ee

Put $L=D$, and modify $A$ and $B$ so that
$A>\frac{\lambda(f)}{\rho(f)}$ and $DB>AL>1$. By Claim 1, for
sufficiently large $n$, there exists $t_n$ such that

\be R_n^{A}\leq t_n\leq R_n^{AL}\textrm{ and
}m(t_n,f)>M(R_n^{A},f)^{D}.\label{3.08}\ee From Claim 2 we see
that when $n$ is sufficiently large
$$ M(R_n^{A},f)>M(R_n,f)^{B}.$$
This and (\ref{3.08}) yield that there exists a positive integer
$n_0$ such that, for $n\geq n_0$, \be m(t_n,f)>M(R_n,f)^{DB}\geq
R_{n+1}^{AL}.\label{3.09}\ee By combining (\ref{3.08}),
(\ref{3.09}), Lemma \ref{lem2.1} and the proof of theorem 3, we have
 proved Theorem 4.  \qed

\vspace{.4cm} \noindent\textbf{\large Acknowledgment: }We would like
to thank Xinhou Hua, Chung-Chun Yang and Jian-Hua Zheng for their
helpful conversations on the subject of this paper. Also we want to
thank the referee for his/her carefully reading over the manuscript,
with helpful suggestions and remarks.

\

\

\

\

\begin{tabular} {ll}
Xiaoling WANG& Wang ZHOU   \\
Department of Applied Mathematics&Department of Statistics and Applied Probability   \\
Nanjing University of Finance and Economics&National University of Singapore    \\
Nanjing 210046, Jiangsu, China & Singapore, 117546                    \\
Email:{\it wangxiaoling@vip.163.com} & Email: {\it stazw@nus.edu.sg} \\

\end{tabular}

\end{document}